\documentclass[12pt]{article}

\usepackage{graphicx,tikz,color,url}
\usepackage{amsmath,amssymb,latexsym}
\usepackage{pgfplots}
\usepackage{mathrsfs}
\usetikzlibrary{arrows}

\newtheorem{theorem}{Theorem}[section]

\newtheorem{lemma}[theorem]{Lemma}
\newtheorem{proposition}[theorem]{Proposition}
\newtheorem{corollary}[theorem]{Corollary}

\newcommand{\s}{\color{blue}}

\newcommand{\proof}{\noindent{\bf Proof.\ }}
\newcommand{\qed}{\hfill $\square$ \bigskip}
\newcommand{\cp}{\,\square\,}
\newcommand{\gp}{{\rm gp}}
\newcommand{\ip}{{\rm ip}}

\newcommand{\diam}{{\rm diam}}

\begin{document}

\title{The general position number under vertex and edge removal}

\author{Pakanun Dokyeesun$^{a}$\thanks{Email: \texttt{papakanun@gmail.com}}
	\and Sandi Klav\v zar $^{a,b,c}$\thanks{Email: \texttt{sandi.klavzar@fmf.uni-lj.si}}
	\and Jing Tian $^{d,b}$\thanks{Email: \texttt{jingtian526@126.com}}
}
\maketitle

\begin{center}
	$^a$ Faculty of Mathematics and Physics, University of Ljubljana, Slovenia\\
	\medskip

	$^b$ Institute of Mathematics, Physics and Mechanics, Ljubljana, Slovenia\\
	\medskip
	
	$^c$ Faculty of Natural Sciences and Mathematics, University of Maribor, Slovenia\\
	\medskip
	
	$^d$ School of Science, Zhejiang University of Science and Technology, Hangzhou, Zhejiang 310023, PR China\\
	\medskip

\end{center}

\begin{abstract}
Let ${\rm gp}(G)$ be the general position number of a graph $G$. It is proved that ${\rm gp}(G-x)\leq 2{\rm gp}(G)$ holds for any vertex $x$ of a connected graph $G$ and that if $x$ lies in some ${\rm gp}$-set of $G$, then ${\rm gp}(G) - 1 \le {\rm gp}(G-x)$. Constructions are given which show that ${\rm gp}(G-x)$ can be much larger than ${\rm gp}(G)$ also when $G-x$ is connected. For diameter $2$ graphs it is proved that ${\rm gp}(G-x) \le {\rm gp}(G)$, and that ${\rm gp}(G-x) \ge {\rm gp}(G) - 1$ when the diameter of $G-x$ remains $2$. It is demonstrated that ${\rm gp}(G)/2\le {\rm gp}(G-e)\leq 2{\rm gp}(G)$ holds for any edge $e$ of a graph $G$. For diameter $2$ graphs these results can be improved to ${\rm gp}(G)-1\le {\rm gp}(G-e)\leq\ {\rm gp}(G) + 1$. All these bounds are proved to be sharp.
\end{abstract}

\noindent
{\bf Keywords:} general position set; vertex-deleted subgraph; edge-deleted subgraph; graph diameter \\

\noindent
{\bf AMS Subj.\ Class.\ (2020)}:  05C12, 05C69

\section{Introduction}

Local operations are valuable in graph theory for understanding and analyzing the properties of graphs and refer to operations that affect only a small part of a graph, rather than the whole structure. These operations include vertex/edge removal/addition and edge subdivision or contraction, and often lead to respective criticality concepts, cf.~\cite{almalki-2024, cao-2022, huo-2022, jakovac-2023, mubayi-2010, zamime-2024}.

The general position problem was introduced to graph theory in~\cite{chandran-2016, manuel-2018a} and proved to be NP-hard in~\cite{manuel-2018a}. In~\cite{anand-2019}, the structure of general position sets was clarified. Afterwards, the problem of determining the general position number of a graph received wide attention, cf.~\cite{ghorbani-2021, Patkos-2019, ThaChaTuiThoSteErs-2024, thomas-2020, tian-2025, tian-2024+, yao-2022}. This is particularly the case for graph products. The general position number of the Cartesian product of paths~\cite{Klavzar-Rus-2021}, of paths and cycles~\cite{Korze-Vesel-2023}, and of two trees~\cite{tian-2021} were determined, while in~\cite{klavzar-2019} the general position number of strong product graphs was investigated. The concept has been modified and/or generalized into several directions. Let us point here to $d$-position sets~\cite{klavzar-2021}, Steiner general position problem~\cite{klavzar-2021a},
edge general position problem~\cite{Klavzar-Tan-2023, manuel-2022, tian-2024}, monophonic general position problem~\cite{thomas-2024+}, and general position polynomials~\cite{Irsic-2023}.

In this paper, we focus on how much the general position number of a graph can be affected by removing a vertex or by removing an edge. In the next section we give definitions and recall some known results that will be needed later on.
In Section~\ref{sec:general}, we prove that $\gp(G-x)\leq 2\gp(G)$ holds for any vertex $x$ of a graph $G$.
Then we demonstrate that $\gp(G-x)$ cannot be bounded from below by a function of $\gp(G)$. On the other hand, if $x$ lies in some $\gp$-set of $G$, then we prove that $\gp(G) - 1 \le \gp(G-x)$.
In Section~\ref{sec:constructions}, we give two constructions which show that $\gp(G-x)$ can be much larger than $\gp(G)$ also when $G-x$ is connected.
In Section~\ref{sec:vertex-diameter-2}, we focus on the vertex removing operation in diameter $2$ graphs.
We show that if $\diam(G) = 2$, then $\gp(G-x) \le \gp(G)$ and prove that $\gp(G) - 1 \le \gp(G-x) \le \gp(G)$ when the diameter of $G-x$ remains $2$.
In Section~\ref{sec:edge remove}, we prove that $\gp(G)/2\le \gp(G-e)\leq 2\gp(G)$ holds for any edge $e$ of a graph $G$. For diameter $2$ graphs $G$ we sharpen the bound by proving that $\gp(G)-1\le \gp(G-e)\leq\ \gp(G) + 1$. All the above bounds are along the way shown to be sharp. We conclude the paper by listing several directions for future investigations. 

\section{Preliminaries}

Unless stated otherwise, the graphs $G=(V(G),E(G))$ considered in this paper are simple and connected. For a positive integer $k$, we use $[k]$ to represent the set $\{1,\ldots, k\}$.
The {\em degree} of a vertex $u$ is the number of vertices adjacent to $u$ in $G$. Vertices of degree one are called {\em leaves}. The number of leaves of $G$ is denoted by $\ell(G)$. If $S\subseteq V(G)$, the subgraph of $G$ induced by $S$ is denoted by $G[S]$. In particular, $G-v$ denotes $G[V(G)\setminus \{v\}]$.
A vertex subset $S$ is an {\em independent set} of $G$ if $G[S]$ is an edgeless graph. The {\em independence number} of $G$, denoted by $\alpha(G)$, is the maximum cardinality of an independent set in $G$.

The {\em distance} $d_G(u,v)$ between vertices $u$ and $v$ of $G$ is the number of edges on a shortest $u,v$-path. A shortest path of $G$ is also called a {\em geodesic} of $G$. The {\em diameter} of $G$ is the maximum distance between pairs of vertices of $G$ and is denoted by $\diam(G)$. A subgraph $H$ of $G$ is {\em isometric} if for each pair of vertices $u,v\in V(H)$ we have $d_H(u,v) = d_G(u,v)$.
The {\em interval} between vertices $u$ and $v$ is
$$I_G[u,v] = \{w:\ d_G(u,v) = d_G(u,w) + d_G(w,v)\}\,.$$

A set $X\subseteq V(G)$ is a {\em general position set} of $G$ if for each pair $u,v \in X$ and any shortest $u,v$-path $P$ we have $V(P)\cap X = \{u,v\}$. The cardinality of a largest general position set of $G$ is the {\em general position number} of $G$ denoted by $\gp(G)$ and referred to as the {\em gp-number of $G$}. A general position set $X$ of cardinality $\gp(G)$ is referred to as a {\em gp-set} of $G$.  

Subgraphs $H_1, \ldots, H_k$ of a graph $G$ form an {\em isometric cover} of $G$ if each $H_i, \, i \in [k]$, is isometric in $G$, and $\bigcup_{i=1}^{k} V(H_i)=V(G)$.

\begin{theorem} {\rm \cite[Theorem 3.1]{manuel-2018a}}
\label{thm:isometric cover lemma}
If $\{H_1, $\ldots$, H_k\}$ is an isometric cover of $G$, then
$$\gp(G) \le \sum_{i=1}^{k}\gp(H_i)\,.$$	
\end{theorem}

The \emph{isometric-path number} of a graph $G$, denoted by $\ip(G)$, is the minimum number of isometric paths required to cover the vertices of $G$.

\begin{proposition}
{\rm\cite[Corollary 3.2]{manuel-2018a}}
\label{prop:ip-ic}
If $G$ is a graph, then $\gp(G) \le 2 \ip(G)$.
\end{proposition}

The following result will be used several times, either implicitly or explicitly.

\begin{proposition}{\rm\cite[Corollary 3.7]{manuel-2018a}}
\label{prop:leaf}
If $T$ is a tree, then $\gp(T)= \ell(T)$.
\end{proposition}

The {\em fan graph} $F_n$, $n\ge 3$, is obtained by taking the join of the path graph $P_n$ and the graph $P_1$. Equivalently, a fan graph is obtained from a wheel graph by removing an edge of it between two degree $3$ vertices, cf.~\cite{selig-2023}.

\begin{proposition}{\rm \cite[Corollary 2.9]{tian-2023}}
\label{prop:fan}
If  $n\ge 4$, then $\gp(F_n) = \lceil \frac{2(n+1)}{3}\rceil$.
\end{proposition}

The final known result we recall describes general position sets in an arbitrary graph. To state it, some more definitions are required. If ${\cal P} = \{S_1, \ldots, S_t\}$ is a partition of $S\subseteq V(G)$, then ${\cal P}$ is \emph{distance-constant} if for any $i,j\in [t]$, $i\ne j$,
there exists a constant $p_{ij}$, such that $d_G(x,y) = p_{ij}$ for every $x\in S_i$, $y\in S_j$. If so, we set $d_G(S_i,S_j) = p_{ij}$. A distance-constant partition ${\cal P}$ is {\em in-transitive} if $p_{ik} \ne p_{ij} + p_{jk}$ holds for $i,j,k\in [t]$.

\begin{theorem} {\rm \cite[Theorem 3.1]{anand-2019}}
\label{thm:gpsets}
Let $G$ be a graph. Then $S\subseteq V(G)$ is a general position set if and only if the components of $G[S]$ are complete subgraphs, the vertices of which form an in-transitive, distance-constant partition of $S$.
\end{theorem}

\section{General bounds}
\label{sec:general}

In this section we prove that $\gp(G-x)\leq 2\gp(G)$ holds for any vertex $x$ of a graph $G$. Then we demonstrate that $\gp(G-x)$ cannot be bounded from below by a function of $\gp(G)$. On the other hand, if $x$ lies in some gp-set, then $\gp(G) - 1 \le \gp(G-x)$.

\begin{theorem}
\label{thm:upper bound}
If $x$ is a vertex of a graph $G$, then $\gp(G-x)\leq 2\gp(G)$. Moreover, the bound is sharp.
\end{theorem}

\proof
Let $R$ be an arbitrary $\gp$-set of $G-x$. Then clearly $x\not\in R$. We partition $R$ into two sets $R_1$ and $R_2$ as follows. For every $u\in R$, we put each vertex from $(I_G[u,x]\setminus \{u\})\cap R$ into $R_2$. This defines $R_2$, and then $R_1 = R\setminus R_2$. We claim that $R_1$ and $R_2$ are general position sets of $G$. 

Suppose first that $R_1$ is not a general position set of $G$. Then there exist vertices $u, v, w\in R_1$ and a shortest $u,w$-path $P$ in $G$ that passes through $v$. Since $R_1\subseteq R$ and $R$ is a general position set of $G-x$, the path $P$ must contain the vertex $x$. We may without loss of generality assume that $x$ lies in the $u,v$-subpath of $P$. Then the $w,x$-subpath of $P$ is a shortest $w,x$-path in $G$. By definition of $R_2$ we get that $v\in R_2$, a contradiction.  

Suppose second that $R_2$ is not a general position set of $G$. Hence there exist vertices $u, v, w\in R_2$ and a shortest $u,w$-path $P$ in $G$ that passes through $v$. Just as in the above paragraph, the path $P$ must contain the vertex $x$ and we may assume that $x$ lies in the $u,v$-subpath of $P$. Since $w\in R_2$, there exists a vertex $w'\in R$, such that $w$ lies on a shortest $w',x$-path $Q$ in $G$. Since $P$ and $Q$ are shortest paths, we infer that the $w,x$-subpath of $Q$ and the $w,x$-subpath of $P$ are of the same length. But this in turn implies that the vertices $w', w, v$ from $R$ lie on a common shortest path in $G-x$, a contradiction.  

We have thus proved that $R_1$ and $R_2$ are general position sets of $G$. Therefore, 
$$\gp(G)\geq \max\{|R_1|, |R_2|\} \geq \frac{1}{2} |R| = \frac{1}{2}\gp(G-x)\,,$$
hence $\gp(G-x) \le 2\gp(G)$.

To show that the bound is sharp, consider the subdivided graph $S(K_{1,n})$, $n\ge 2$, of the star $K_{1,n}$, that is, the graph obtained from $K_{1,n}$ by subdividing each of its edges once. Let $x$ be the vertex of degree $n$ of $S(K_{1,n})$. Since $S(K_{1,n})-x \cong nK_2$ and having Proposition~\ref{prop:leaf} in mind, we can conclude that $\gp(S(K_{1,n})-x) = 2n = 2\gp(S(K_{1,n}))$.
\qed

There is no general lower bound on $\gp(G-x)$ in terms of $\gp(G)$. To demonstrate it, consider the fan graphs $F_n$, $n\ge 3$. By Proposition~\ref{prop:fan} we have $\gp(F_n) = \lceil \frac{2(n+1)}{3}\rceil$. Since clearly $\gp(F_n - x) = 2$, where $x$ is the vertex of $F_n$ of degree $n$, we see that $\gp(G-x)$ can indeed be arbitrarily smaller than $\gp(G)$. The next result leads to many additional such examples.

\begin{proposition}
\label{prop:much-smaller}
Let $S$ be an independent set of a graph $H$ with $|S| = \alpha(H)$. If $G$ is the graph obtained from the disjoint union of $H$ and a vertex $x$ by joining $x$ to each vertex of $S$, then $\gp(G) \ge \alpha(H)$.
\end{proposition}

\proof
In $G$, the set $S$ is an independent set of vertices that are pairwise at distance $2$. Hence $S$ is a general position set of $H$ and the conclusion follows.
\qed

In Proposition~\ref{prop:much-smaller} we have $H \cong G - x$, hence $\gp(G-x) = \gp(H)$. Thus, if $\gp(H)$ is much smaller than $\alpha(H)$, then $\gp(G-x)$ is much smaller than $\gp(G)$. For instance, such graphs are grids $P_n\cp P_m$, $n,m\ge 3$, for which we know that $\gp(P_n\cp P_m) = 4$~\cite[Corollary 3.2]{manuel-2018b}.

On the other hand, under some additional assumption, $\gp(G-x)$ can be bounded from below with $\gp(G)$ as follows.

\begin{proposition}\label{pro: in_gp-set}
Let $x$ be a vertex of a graph $G$. If $x$ lies in some gp-set of $G$, then $\gp(G) -1 \le \gp(G-x)$.
\end{proposition}

\proof
Let $S$ be a gp-set of $G$ and $x\in S$.
Suppose that $S \setminus \{x\}$ is not a general position set of $G-x$. Then there are three distinct vertices  $u, v, w \in S\setminus \{x\}$ lying on a shortest path in $G-x$. Without loss of generality, assume that $v$ lies on a $u,w$-geodesic in $G-x$.
Since $S$ is a gp-set of $G$, but $S \setminus \{x\}$ is not a general position set of $G-x$, $x$ must lie on a $u,w$-geodesic in $G$. This contradicts our assumption. Hence, $S \setminus \{x\}$ is a general position set of $G-x$. It concludes that $ \gp(G-x) \ge |S|-1 =\gp(G) -1.$
\qed

\section{Two constructions}
\label{sec:constructions}

In Theorem~\ref{thm:upper bound} we have proved that the bound $\gp(G-x)\leq 2\gp(G)$ is sharp. However, sharpness examples were such that $G-x$ is not connected. In this section we give two constructions which show that $\gp(G-x)$ can be much larger than $\gp(G)$ also when $G-x$ is connected.

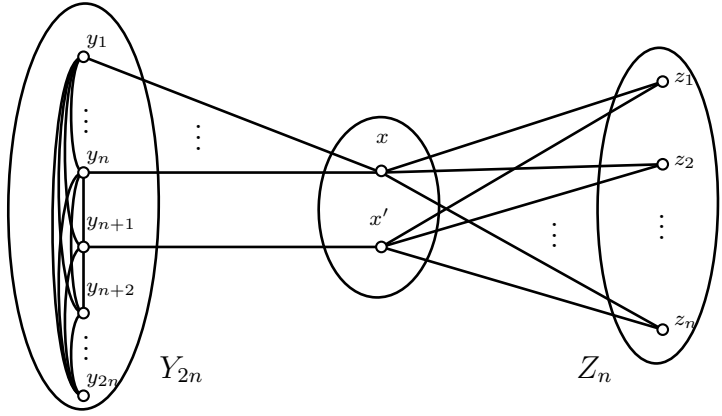
\begin{figure}[ht!]
	\begin{center}
\begin{tikzpicture}
[line cap=round,line join=round,>=triangle 45,x=1.1cm,y=1.1cm, style=thick]
	\clip(3.22,-0.22) rectangle (13.42,5.26);
	\draw [line width=1.0pt] (12.,4.)-- (8.6,2.9);
	\draw [line width=1.0pt] (8.6,2.9)-- (12.,3.);
	\draw [line width=1.0pt] (8.6,2.9)-- (12.,1.);
	\draw [line width=1.0pt] (8.6,2.0)-- (12.,4.);
	\draw [line width=1.0pt] (8.6,2.0)-- (12.,3.);
	\draw [line width=1.0pt] (8.6,2.0)-- (12.,1.);
	\draw [rotate around={89.48615028364591:(5.1,2.49)},line width=1.0pt] (5.1,2.49) ellipse (2.8cm and 1cm);
	\draw [line width=1.0pt] (5.,2.9)-- (8.6,2.9);
	\draw [line width=1.0pt] (5.,4.3)-- (8.6,2.9);
	\draw [line width=1.0pt] (5.,2.)-- (8.6,2.0);
	\draw [line width=1.0pt,color=black]  (5.,4.3)..controls(4.8,4.3)and (4.8,2.9)..(5.,2.9);
	\draw [line width=1.0pt,color=black]  (5.,4.3)..controls(4.7,4.3)and (4.7,2)..(5.,2);
	\draw [line width=1.0pt,color=black]  (5.,4.3)..controls(4.6,4.3)and (4.6,1.2)..(5.,1.2);
	\draw [line width=1.0pt,color=black]  (5.,4.3)..controls(4.5,4.3)and (4.5,0.2)..(5.,0.2);
	\draw [line width=1.0pt,color=black]  (5.,2.9)--(5.,2);
	\draw [line width=1.0pt,color=black]  (5.,2.9)..controls(4.8,2.9)and (4.8,1.2)..(5.,1.2);
	\draw [line width=1.0pt,color=black]  (5.,2.9)..controls(4.6,2.9)and (4.6,0.2)..(5.,0.2);
	\draw [line width=1.0pt,color=black]  (5.,2)--(5.,1.2);
	\draw [line width=1.0pt,color=black]  (5.,2.)..controls(4.7,2.)and (4.7,0.2)..(5.,0.2);
		\draw [line width=1.0pt,color=black]  (5.,1.2)..controls(4.8,1.2)and (4.8,0.2)..(5.,0.2);
	\draw (4.84,4) node[anchor=north west] {$\vdots$};
	\draw (6.2,3.8) node[anchor=north west] {$\vdots$};
	\draw (4.84,1.25) node[anchor=north west] {$\vdots$};
	\draw [rotate around={89.38394009160086:(12.1,2.5)},line width=1.0pt] (12.1,2.5) ellipse (2.2cm and 0.8cm);
	\draw (4.9,4.8) node[anchor=north west]{\small$y_1$};
	\draw (4.9,3.4) node[anchor=north west] {\small$y_n$};
	\draw (4.9,2.6) node[anchor=north west] {\small$y_{n+1}$};
	\draw (4.9,1.7) node[anchor=north west] {\small$y_{n+2}$};
	\draw (4.9,0.7) node[anchor=north west] {\small$y_{2n}$};
	\draw (8.4,3.5) node[anchor=north west] {\small$x$};
	\draw (8.32,2.6) node[anchor=north west] {\small$x'$};
	\draw (12,4.24) node[anchor=north west] {\small$z_1$};
	\draw (12,3.24) node[anchor=north west] {\small$z_2$};
	\draw (12,1.3) node[anchor=north west] {\small$z_n$};
	\draw (5.78,0.82) node[anchor=north west] {$Y_{2n}$};
	\draw (10.82,0.82) node[anchor=north west] {$Z_n$};
	\draw [rotate around={88.0474909506011:(8.57,2.48)},line width=1.0pt] (8.57,2.48) ellipse (1.2cm and 0.8cm);
	\draw (11.8,2.62) node[anchor=north west] {$\vdots$};
	\draw (10.5,2.62) node[anchor=north west] {$\vdots$};
		\draw [fill=white] (12.,4.) circle (2pt);
		\draw [fill=white] (8.6,2.92) circle (2pt);
		\draw [fill=white] (12.,3.) circle (2pt);
		\draw [fill=white] (12.,1.) circle (2pt);
		\draw [fill=white] (8.6,2) circle (2pt);
		\draw [fill=white] (5,4.3) circle (2pt);
		\draw [fill=white] (5,2.9) circle (2pt);
		\draw [fill=white] (5,2.) circle (2pt);
		\draw [fill=white] (5,1.2) circle (2pt);
		\draw [fill=white] (5,0.2) circle (2pt);

\end{tikzpicture}
	\caption{Graph $H_n$.}
\label{fig:H_n}
\end{center}	
\end{figure}

In the first construction let $H_n$, $n\ge 3$, be the graph defined as follows. Its vertex set is $Y_{2n} \cup \{x,x'\} \cup Z_n$, where $Y_{2n} = \{y_1, \ldots, y_{2n}\}$ and $Z_n = \{z_1, \ldots, z_n\}$. The vertices of $Y_{2n}$ induce a complete subgraph $K_{2n}$, the vertices of $\{x,x'\} \cup Z_n$ induce a complete bipartite graph $K_{2,n}$ with the corresponding bipartition, the vertex $x$ is adjacent to vertices $y_1, \ldots, y_n$, and the vertex $x'$ is adjacent to $y_{n+1}$. See Fig.~\ref{fig:H_n}.

\begin{proposition}
\label{prop:Hn}
If $n\ge 3$, then $\gp(H_n) = 2n + 1$ and $\gp(H_n-x) = 3n-1$.
\end{proposition}

\proof
Set $Y_{j}=\{y_1,\ldots,y_{j}\}$ for $j\in \{n,n+1, 2n\}$.  Using Theorem~\ref{thm:gpsets}, we see that $Y_{n+1}\cup Z_n$ is a general position set of $H_n$. It follows that $\gp(H_n) \geq |Y_{n+1}\cup Z_n| = 2n + 1$.

To prove the upper bound on $\gp(H_n)$, suppose on the contrary that $\gp(H_n) \geq 2n + 2$ and let $R$ be a gp-set of $H_n$. We claim first that $|Y_n\cap R|\geq 1$ and $|Z_n\cap R|\geq 1$. Indeed, if $|Y_n\cap R| = 0$, then $R = \{y_{n+1}, \ldots, y_{2n}\} \cup \{x,x'\} \cup Z_n$, but this is clearly not a general position set. Similarly, if $|Z_n\cap R| = 0$, then $R = Y_{2n} \cup \{x,x'\}$, which is also not a general position set as $y_{2n}$, $y_{n+1}$, and $x'$ lie on a common shortest path. Hence the claim. Since $|Y_n\cap R|\geq 1$ and $|Z_n\cap R|\geq 1$, we get that $x\not\in R$ and $(Y_{2n}\setminus Y_{n+1})\cap R = \emptyset$. It follows that $R = Y_{n+1}\cup \{x'\} \cup Z_n$. But then $z_n$, $x'$, and $y_{n+1}$ are on a shortest path. This final contradiction proves that $\gp(H_n) \le 2n + 1$. We have thus shown that  $\gp(H_n) = 2n + 1$.

Consider now $H_n-x$ and note that $\diam(H_n-x)=\diam(H_n)=3$. In~\cite[Theorem 2.4]{chandran-2016} it was proved that if $G$ is a graph, then $\gp(G)\leq n(G)-\diam(G)+1$. Hence $\gp(H_n-x)\leq 3n-1$. Invoking Theorem~\ref{thm:gpsets} again, we infer that $(Y_{2n}\setminus\{y_{n+1}\})\cup Z_n$ is a general position set of $H_n-x$ which in turn implies that $\gp(H_n-x)=3n-1$.
\qed

\begin{figure}[ht!]
\begin{center}
\begin{tikzpicture}
[line cap=round,line join=round,>=triangle 45,x=1.0cm,y=1.0cm,style=thick]
	\clip(-3.56,-1.) rectangle (7.26,5.2);
	\draw [line width=1.0pt] (-2.,4.)-- (0.16,3);
	\draw [line width=1.0pt] (-2.,3.)-- (0.16,3);
	\draw [line width=1.0pt] (0.16,3)-- (-2.,1.);
	\draw [line width=1.0pt] (0.16,3)-- (2.,4.);
	\draw [line width=1.0pt] (0.16,3)-- (2.,3.);
	\draw [line width=1.0pt] (0.16,3)-- (2.,1.);
	\draw [line width=1.0pt] (2.,4.)-- (5.5,4.);
	\draw [line width=1.0pt] (2.,1.)-- (5.5,3.);
	
	\draw [line width=1.0pt,color=black]  (5.5,4)--(5.5,3);
	\draw [line width=1.0pt,color=black]  (5.5,4)..controls(5.8,4)and (5.8,1.5)..(5.5,1.5);
	\draw [line width=1.0pt,color=black]  (5.5,4)..controls(5.9,4)and (5.9,0.5)..(5.5,0.5);
	\draw [line width=1.0pt,color=black]  (5.5,1.5)--(5.5,0.5);
	\draw [line width=1.0pt,color=black]  (5.5,3)..controls(5.65,3)and (5.65,1.5)..(5.5,1.5);
	\draw [line width=1.0pt,color=black]  (5.5,3)..controls(5.8,3)and (5.8,0.5)..(5.5,0.5);
	\draw [rotate around={88.7939512207795:(-2.,2.5)},line width=1.0pt] (-2.,2.5) ellipse (2.0cm and 0.9cm);
	\draw [rotate around={89.39371808448271:(1.96,2.47)},line width=1.0pt] (1.96,2.47) ellipse (2cm and 0.9cm);
	\draw (-0.12,2.86) node[anchor=north west] {\small $w$};
	\draw (-2.5,4) node[anchor=north west] {\small $x_1$};
	\draw (-2.5,3) node[anchor=north west] {\small $x_2$};
	\draw (-2.5,1.5) node[anchor=north west] {\small $x_k$};
	\draw (1.85,4) node[anchor=north west] {\small $y_1$};
	\draw (1.85,3) node[anchor=north west] {\small $y_2$};
	\draw (1.85,1.7) node[anchor=north west] {\small $y_k$};
	\draw (5,4) node[anchor=north west] {\small $z_1$};
	\draw (5,2.9) node[anchor=north west] {\small $z_2$};
	\draw (4.9,1.7) node[anchor=north west] {\small $z_k$};
	\draw (4.95,0.55) node[anchor=north west] {\small$z_{k+1}$};
	\draw (-2.2,0.24) node[anchor=north west] {$X_k$};
	\draw (1.72,0.24) node[anchor=north west] {$Y_k$};
	\draw (4.8,0) node[anchor=north west] {$Z_{k+1}$};
	\draw [rotate around={-88.78112476486815:(5.5,2.4)},line width=1.0pt] (5.5,2.4) ellipse (2.5cm and 1cm);
	\draw [line width=1.0pt] (2.,3.)-- (5.46,3.04);

	\draw (5.3,2.64) node[anchor=north west] {$\vdots$};
	\draw (1.9,2.54) node[anchor=north west] {$\vdots$};
	\draw (-2.16,2.62) node[anchor=north west] {$\vdots$};
	\draw (3.5,3) node[anchor=north west] {$\vdots$};

		\draw [fill=white] (-2.,4.) circle (2pt);
		\draw [fill=white] (-2.,3.) circle (2pt);
		\draw [fill=white] (-2.,1.) circle (2pt);
		\draw [fill=white] (0.16,3) circle (2pt);
		\draw [fill=white] (2.,4.) circle (2pt);
		\draw [fill=white] (2.,3.) circle (2pt);
		\draw [fill=white] (2.,1.) circle (2pt);
		\draw [fill=white] (5.5,4.) circle (2pt);
		\draw [fill=white] (5.5,3.) circle (2pt);
		\draw [fill=white] (5.5,1.5) circle (2pt);
		\draw [fill=white] (5.5,0.5) circle (2pt);

\end{tikzpicture}
	\caption{Graph $G_k$.}
	\label{fig:G_k}
\end{center}	
\end{figure}
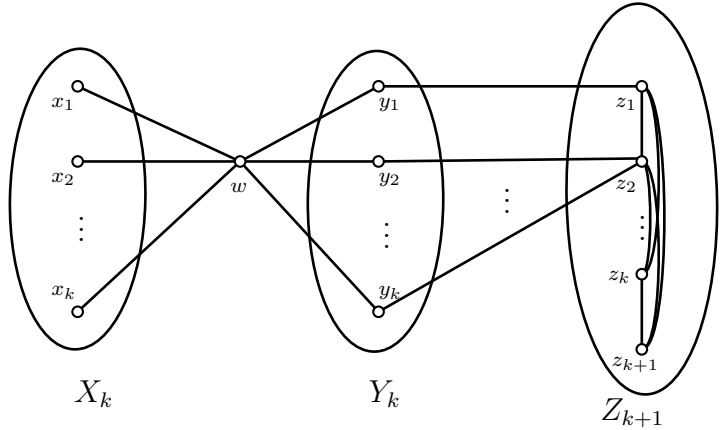

We next give another family of graphs in which the general position number increases arbitrarily by removing a vertex. If $k \ge 2$, then let the graph $G_k$ be constructed as follows. Let $V(G_k) = X_k \cup Y_k \cup Z_{k+1} \cup \{w\}$, where $X_k = \{x_1, \ldots, x_k\}$, $Y_k = \{y_1, \ldots, y_k\}$, and $Z_{k+1} = \{z_1, \ldots, z_{k+1}\}$. The vertex $w$ is adjacent to every vertex of $X_k\cup Y_k$, the vertices of $Z_{k+1}$ induce a complete graph $K_{k+1}$, $z_2$ is adjacent to $y_2, \ldots, y_k$, and $z_1$ is adjacent to $y_1$, see Fig.~\ref{fig:G_k}. Note that $\diam(G_k) = 4$, and that $X_k$ and $Y_k$ are independent sets of vertices.

\begin{proposition}
\label{prop:Gk}
If $k \ge 2$, then $\gp(G_k) = 2k$ and $\gp(G_k - z_2)=3k-2$.
\end{proposition}

\proof
Let $P_1$ be the path induced by the vertices $x_1,w, y_1$, and $z_1$. Let $P_i=x_iwy_iz_2z_{i+1}$ be a path of $G_k$, where $2\leq i\leq k$.
Set $\Psi=\{P_i:\ i\in [k]\}$. Then it follows that $|\Psi|=k$. It is observed that $\Psi$ is a set of isometric paths of $G_k$.
By Proposition~\ref{prop:ip-ic}, $\gp(G_k) \le 2k$.
To show $\gp(G_k) \ge 2k$, note that $G[X_k \cup Y_k\cup \{w\}] \cong K_{1,2k}$, hence Proposition~\ref{prop:leaf} implies that $\gp(G[X_k \cup Y_k\cup \{w\}]) = 2k$. Since $G[X_k \cup Y_k\cup \{w\}]$ is an isometric subgraph of $G_k$, we conclude that $\gp(G_k) \ge 2k$ and thus $\gp(G_k) = 2k$.

Consider $G_k -z_2$. It is straightforward to check that $S = V(G_k-z_2) \setminus \{w, y_1, z_1\}$ is a general position set of
$G_k -z_2$. For instance, any shortest path between a vertex $x_i\in X_k$ and any other vertex from $S$ avoids other vertices of $S$ because $d_{G_k -z_2}(x_i, x_j) = 2$, $d_{G_k -z_2}(x_i, y_j) = 2$, and $d_{G_k -z_2}(x_i, z_j) = 4$ for $j\ge 3$. In the latter case, a shortest $x_i,z_j$-path, induced by the vertices $x_i, w, y_1, z_1$, and $z_j$, is unique. Thus $\gp(G_k -z_2) \ge 3k -2$. On the other hand, suppose that $S$ is a general position set of $G_k -z_2$ of size $3k-1$. Because $n(G_k -z_2) = 3k+1$, we have $S\cap \{w, y_1, z_1\}\ne \emptyset$. But then we find a vertex from $X_k$, a vertex from $\{w, y_1, z_1\}$, and a vertex from $Z_{k+1}\setminus\{z_1, z_2\}$, such that these three vertices from $S$ lie on a common shortest path in $G_k -z_2$. As a consequence, we conclude that $\gp (G_k - z_2) \le 3k-2$ and thus we have $\gp (G_k - z_2) = 3k-2$.
\qed

\section{Vertex removing in diameter $2$ graphs}
\label{sec:vertex-diameter-2}

Note that for the graphs $H_n$ from Proposition~\ref{prop:Hn} we have $\diam(H_n) = 3$ as well as $\diam(H_n - x) = 3$. Also, for the graphs $G_k$ from Proposition~\ref{prop:Gk} we have $\diam(G_k) = 4$ as well as $\diam(G_k - z_2) = 4$. In both cases we have seen that removing a vertex increases the general position number arbitrarily. In this section we therefore focus on the vertex removing operation in diameter $2$ graphs. We first show that in this case Theorem~\ref{thm:upper bound} can be sharpened to $\gp(G-x) \le \gp(G)$. Second, we prove that if the diameter of $G-x$ remains $2$, then $\gp(G) - 1\le \gp(G-x) \le \gp(G)$. Before presenting these results, we consider some examples.

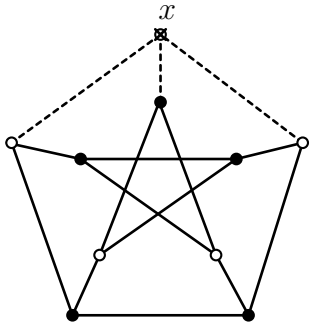
\begin{figure}[ht!]
	\begin{center}
		\usetikzlibrary{shapes.misc}
		
		\tikzset{cross/.style={cross out, draw=black, minimum size=4*(#1-\pgflinewidth), inner sep=0pt, outer sep=0pt},
			cross/.default={1pt}}
\begin{tikzpicture}[line cap=round,line join=round,>=triangle 45,x=0.9cm,y=0.9cm, style=thick]
	\draw [line width=1.0pt] (4.2,0.25)-- (6.8,0.25);
	\draw [line width=1.0pt] (6.8,0.25)-- (7.6,2.8);
	\draw [line width=1.0pt,dash pattern=on 2pt off 2pt] (7.6,2.8)-- (5.5,4.4);
	\draw [line width=1.0pt,dash pattern=on 2pt off 2pt] (5.5,4.4)-- (3.3,2.8);
	\draw [line width=1.0pt] (3.3,2.8)-- (4.2,0.25);
	\draw [line width=1.0pt,dash pattern=on 2pt off 2pt] (5.5,3.4)-- (5.5,4.4);
	\draw [line width=1.0pt] (6.62,2.56)-- (7.6,2.8);
	\draw [line width=1.0pt] (6.32,1.14)-- (6.8,0.25);
	\draw [line width=1.0pt] (4.6,1.14)-- (4.2,0.25);
	\draw [line width=1.0pt] (3.23,2.8)-- (4.32,2.56);
	\draw [line width=1.0pt] (4.6,1.14)-- (5.48,3.4);
	\draw [line width=1.0pt] (5.48,3.4)-- (6.32,1.14);
	\draw [line width=1.0pt] (6.32,1.14)-- (4.32,2.56);
	\draw [line width=1.0pt] (4.32,2.56)-- (6.62,2.56);
	\draw [line width=1.0pt] (6.62,2.56)-- (4.6,1.14);
	\draw (5.3,5) node[anchor=north west] {$x$};

		\draw [fill=black] (4.2,0.25) circle (2pt);
		\draw [fill=black] (6.8,0.25) circle (2pt);
		\draw [fill=white] (7.6,2.8) circle (2pt);
		\draw [fill=white] (5.5,4.4) circle (2pt);
		\draw [fill=white] (3.3,2.8) circle (2pt);
		\draw [fill=black] (5.5,3.4) circle (2pt);
		\draw (5.5,4.4) node[cross=2pt,line width=1.0pt,black]{};
		\draw [fill=black] (6.62,2.56) circle (2pt);
		\draw [fill=white] (6.32,1.14) circle (2pt);
		\draw [fill=white] (4.6,1.14) circle (2pt);
		\draw [fill=black] (4.32,2.56) circle (2pt);

\end{tikzpicture}
	\caption{A gp-set of $P-x$.}
\label{fig:petersen}
\end{center}	
\end{figure}

The gp-number of a diameter $2$ graph may stay the same after a vertex is removed. Consider for instance complete bipartite graphs $K_{n,m}$, where $2\le n\le m$. Then it is known that $\gp(K_{n,m}) = m$, see~\cite[Proposition 2.2]{chandran-2016}. Hence if $x$ is a vertex of $K_{n,m}$ from a smaller partition set, then $\gp(K_{n,m} - x) = m = \gp(K_{n,m})$. Consider next the Petersen graph $P$. Then $\gp(P) = 6$, see~\cite[page~184]{manuel-2018a}. If $x\in V(P)$, then $\diam(P-x) = 3$. Moreover, by a case analysis we can check that $\gp(P-x) = 5$. See Fig.~\ref{fig:petersen} where a gp-set in $P-x$ is marked with black vertices.

\begin{proposition}
\label{prop:diaG_2}
If $x$ is a vertex of a diameter $2$ graph $G$, then $\gp(G-x) \le \gp(G)$.
\end{proposition}

\proof
To prove the proposition it suffices to show that if $S$ is a gp-set of $G-x$, then $S$ is also a general position set of $G$. Let $u,v,w$ be vertices from $S$ and suppose by way of contradiction that they lie on a shortest path in $G$. As $\diam(G) = 2$, the vertices $u,v,w$ induce an isometric $P_3$ in $G$. But then this path is also isometric in $G-x$, a contradiction.
\qed

\begin{theorem}\label{thm:dia_2}
Let $x$ be a vertex of a diameter $2$ graph $G$. If $\diam(G-x)=2$, then $\gp(G)-1 \le \gp(G-x) \le \gp(G)$. Moreover, the bounds are sharp.
\end{theorem}

\proof
The upper bound follows by Proposition~\ref{prop:diaG_2}. Assume that $\diam(G) =\diam(G-x)=2$, where $x \in V(G)$. Then we can prove along the lines of the proof of Proposition~\ref{prop:diaG_2} that if $S$ is a gp-set of $G$, then $S\setminus \{x\}$ is a general position set of $G-x$. Hence $\gp(G-x) \ge |S\setminus \{x\}| \ge |S| - 1 = \gp(G) -1$.

Let $n_1\ge \cdots \ge n_k\ge 2$, $k\ge 3$, and let $G_{n_1, \ldots, n_k}$ be the graph obtained from the disjoint union of $K_{n_1}, \ldots, K_{n_k}$ by selecting a vertex in each of the complete graphs and identify all of them into a single vertex. Then $\gp(G_{n_1, \ldots, n_k}) = n_1 + \cdots + n_k - k$. If $x$ is an arbitrary vertex of $G_{n_1, \ldots, n_k}$ which is not its maximum degree vertex, then $\gp(G_{n_1, \ldots, n_k}-x) = \gp(G_{n_1, \ldots, n_k}) - 1$. This demonstrates sharpness of the lower bound as these graphs are of diameter 2.

Let next $n_1\ge \cdots \ge n_k\ge 2$, where $k\ge 2$ and $n_1 > k$, and consider the complete multipartite graph $K_{n_1, \ldots, n_k}$. Then $\gp(K_{n_1, \ldots, n_k}) = n_1$ and
if $x$ is an arbitrary vertex which does not lie in the $n_1$-part, then $\gp(K_{n_1, \ldots, n_k}-x) = n_1$, which demonstrates sharpness of the upper bound.
\qed

We note that the bounds from Proposition~\ref{prop:diaG_2} and Theorem~\ref{thm:dia_2} can also be deduced from the fact that a subset of a diameter two graph is in general position if and only if it is an independent union of cliques.

Another family which demonstrates sharpness of the upper bound in Theorem~\ref{thm:dia_2} is the family of strong products $K_n\boxtimes C_m$, where $n\ge 2$ and $m\in \{4,5\}$. By~\cite[Proposition 4.3]{klavzar-2019} we have $\gp(K_n\boxtimes C_4) = 2n$ and $\gp(K_n\boxtimes C_5) = 3n$. Moreover, it is straightforward to check that the general position number does not change after one vertex is removed from these graphs.

Putting together Propositions~\ref{pro: in_gp-set} and~\ref{prop:diaG_2}, and Theorem~\ref{thm:dia_2}, the following conclusion follows.

\begin{corollary}
\label{cor:diam2-all-together}
Let $x$ be a vertex of a diameter $2$ graph $G$. If $\diam(G-x) = 2$, or $x$ lies in some gp-set of $G$, then
$$\gp(G)-1 \le \gp(G-x) \le \gp(G)\,.$$
\end{corollary}

\section{Edge removing in general graphs}
\label{sec:edge remove}

In this section we consider how much the general position number can be affected by removing an edge. Contrary to the vertex removal, we can give general sharp lower and upper bounds. To prove them, we first recall the following well-known sets from metric graph theory. For instance, these sets are of interest in the study of partial cubes~\cite{ovchinnikov-2011}, of distance-balanced graphs~\cite{Jerebic:2008}, of $\ell$-distance-balanced graphs~\cite{Miklavic:2018}, and of the Mostar index~\cite{ali-2021}.

If $e=uv$ is an edge of a graph $G$, then
\begin{align*}
W_{uv} & = \{w\in V(G):\ d_G(u,w) < d_G(v,w) \}, \\
W_{vu} & = \{w\in V(G):\ d_G(v,w) < d_G(u,w) \}, \\
_vW_u & = \{w\in V(G):\ d_G(u,w) = d_G(v,w) \}\,.
\end{align*}
For vertices $x,y\in V(G)$, let further ${\cal P}_G(x,y)$ be the set of all shortest $x,y$-paths in $G$. The following technical lemma about these sets will be crucial for our following arguments.

\begin{lemma}
\label{lem:paths-in-X-u}
Let $e=uv$ be an edge in a graph $G$ and let $x,y\in W_{uv} \cup\, _vW_u$. Then ${\cal P}_{G}(x,y) = {\cal P}_{G-e}(x,y)$. In particular, $d_G(x,y) = d_{G-e}(x,y)$.
\end{lemma}

\proof
Let $P\in {\cal P}_{G}(x,y)$. We claim that $P$ does not contain $e$. Suppose on the contrary that $P$ contains $e$. Assume first that the sequence of the vertices on $P$ is $x, \ldots, v, u, \ldots, y$. Since $d_G(x,u) \le d_G(x,v)$, the path $P$ is not shortest in $G$, a contradiction. Assume second that the sequence of the vertices on $P$ is $x, \ldots, u, v, \ldots, y$. But now the fact that $d_G(y,u) \le d_G(y,v)$ gives another contradiction with the assumption that $P$ is shortest in $G$. We can conclude that $P\in {\cal P}_{G-e}(x,y)$, therefore ${\cal P}_{G}(x,y)  \subseteq {\cal P}_{G-e}(x,y)$.

Let now $P\in {\cal P}_{G-e}(x,y)$. Suppose that  $P\notin {\cal P}_{G}(x,y)$. This means that there exists an $x,y$-path $Q$ in $G$ shorter than $P$. For this to happen, $Q$ must necessarily contain the edge $uv$. But then we can argue analogously as in the first paragraph that $Q$ is not a shortest path. This contradiction implies that $P\in {\cal P}_{G}(x,y)$. Consequently, ${\cal P}_{G-e}(x,y)  \subseteq {\cal P}_{G}(x,y)$, and we are done.
\qed

The main result of this section reads as follows.

\begin{theorem}
\label{thm:edge-removal-upper}
If $e$ is an edge of a graph $G$, then
$$ \frac{\gp(G)}{2}\le \gp(G-e)\leq\ 2\gp(G)\,.$$
Moreover, both bounds are sharp.
\end{theorem}

\proof
Let $e=uv$ and let $X$ be a gp-set of $G$. Setting
\begin{align*}
 X_{uv} & = \{w\in X:\ d_G(u,w) < d_G(v,w) \}, \\
X_{vu} & = \{w\in X:\ d_G(v,w) < d_G(u,w) \}, \\
_vX_u & = \{w\in X:\ d_G(u,w) = d_G(v,w) \},
\end{align*}
we have $X = X_{uv} \cup X_{vu} \cup\, _vX_u$. Let $X_u = X_{uv} \cup\, _vX_u$ and $X_v = X_{vu} \cup\, _vX_u$. We now show that $X_u$ and $X_v$ are general position sets of $G-e$. By symmetry it suffices to prove the claim for $X_u$. Consider any two vertices $x,y\in X_u$ and let $P$ be an arbitrary shortest $x,y$-path in $G-e$. By Lemma~\ref{lem:paths-in-X-u}, the path $P$ is also a shortest $x,y$-path in $G$, hence $V(P) \cap X_u \subseteq \{x,y\}$. It follows that $X_u$ is a general position set of $G-e$ and hence also $X_v$ is such. Therefore
$$\gp(G-e) \ge \max\{ |X_u|, |X_v|\} \ge \frac{|X|}{2} = \frac{\gp(G)}{2}\,.$$
This proves the lower bound.

To prove the upper bound we proceed similarly as above. For this sake let $Y$ be a gp-set of $G-e$ and partition $Y$ into the following subsets:
\begin{align*}
Y_{uv} & = \{w\in Y:\ d_{G-e}(u,w) < d_{G-e}(v,w) \}, \\
Y_{vu} & = \{w\in Y:\ d_{G-e}(v,w) < d_{G-e}(u,w) \}, \\
_vY_u & = \{w\in Y:\ d_{G-e}(u,w) = d_{G-e}(v,w) \}.
\end{align*}
Then, using Lemma~\ref{lem:paths-in-X-u} as above, $Y_u = Y_{uv} \cup\, _vY_u$ and $Y_v = Y_{vu} \cup\, _vY_u$ are general position sets of $G$. Hence
$$\gp(G) \ge \max\{ |Y_u|, |Y_v|\} \ge \frac{|Y|}{2} = \frac{\gp(G-e)}{2}\,,$$
which proves the upper bound.

To prove sharpness of the lower bound, let $k\geq 3$ and let $G_k'$ be the graph constructed as follows. Consider the disjoint union of four copies of $K_{2,k}$, where in two copies an extra edge between its degree $k$ vertices is added (these are the edges $f$ and $f'$ in Fig.~\ref{fig:graph G'}). Then circularly connect these four graphs by three edges and a path of length $8$ as shown in Fig.~\ref{fig:graph G'}.

\begin{figure}[ht!]
\begin{center}
\begin{tikzpicture}[scale=1.0,style=thick,x=1cm,y=1cm]
\def\vr{2pt}

\begin{scope}[xshift=0cm, yshift=0cm] 
\coordinate(x1) at (0,0);
\coordinate(x2) at (1,0.4);
\coordinate(x3) at (2,0);
\coordinate(x4) at (3,0);
\coordinate(x5) at (5,0);
\coordinate(x6) at (6,0);
\coordinate(x7) at (8,0);
\coordinate(x8) at (9,0);
\coordinate(x9) at (10,0.4);
\coordinate(x10) at (11,0);
\coordinate(x11) at (11,2);
\coordinate(x12) at (8.5,2);
\coordinate(x13) at (7,2);
\coordinate(x14) at (5.5,2);
\coordinate(x15) at (4,2);
\coordinate(x16) at (2.5,2);
\coordinate(x17) at (0,2);
\coordinate(x18) at (1,0.75);
\coordinate(x19) at (1,1.5);
\coordinate(x20) at (4,0.75);
\coordinate(x21) at (4,1.5);
\coordinate(x22) at (7,0.75);
\coordinate(x23) at (7,1.5);
\coordinate(x24) at (10,0.75);
\coordinate(x25) at (10,1.5);
\coordinate(w1) at (4,0.4);
\coordinate(w2) at (7,0.4);

\draw (x1) -- (x2) -- (x3) -- (x4);
\draw (x5)-- (x6);
\draw (x7)--(x8) -- (x9) -- (x10) -- (x11) --(x12)-- (x13) --(x14)--(x15) -- (x16) -- (x17) -- (x1);
\draw (x1) -- (x18) -- (x3) -- (x19)--(x1);
\draw (x4) -- (x20) -- (x5) -- (x21)--(x4);
\draw (x6) -- (x22) -- (x7) -- (x23)--(x6);
\draw (x8) -- (x24) -- (x10) -- (x25)--(x8);
\draw (1,1.2) node{\vdots};
\draw (4,1.2) node{\vdots};
\draw (7,1.2) node{\vdots};
\draw (10,1.2) node{\vdots};
\draw (x4) -- (w1) -- (x5);
\draw (x6) -- (w2) -- (x7);
\draw (x4) .. controls (3.5,-0.5) and (4.5,-0.5) .. (x5);
\draw (x6) .. controls (6.5,-0.5) and (7.5,-0.5) .. (x7);

\draw(w1)[fill=white] circle(\vr);
\draw(w2)[fill=white] circle(\vr);
\draw(x1)[fill=white] circle(\vr);
\draw(x2)[fill=white] circle(\vr);
\draw(x3)[fill=white] circle(\vr);
\draw(x4)[fill=white] circle(\vr);
\draw(x5)[fill=white] circle(\vr);
\draw(x6)[fill=white] circle(\vr);
\draw(x7)[fill=white] circle(\vr);
\draw(x8)[fill=white] circle(\vr);
\draw(x9)[fill=white] circle(\vr);
\draw(x10)[fill=white] circle(\vr);
\draw(x11)[fill=white] circle(\vr);
\draw(x12)[fill=white] circle(\vr);
\draw(x13)[fill=white] circle(\vr);
\draw(x14)[fill=white] circle(\vr);
\draw(x15)[fill=white] circle(\vr);
\draw(x16)[fill=white] circle(\vr);
\draw(x17)[fill=white] circle(\vr);
\draw(x18)[fill=white] circle(\vr);
\draw(x19)[fill=white] circle(\vr);
\draw(x20)[fill=white] circle(\vr);
\draw(x21)[fill=white] circle(\vr);
\draw(x22)[fill=white] circle(\vr);
\draw(x23)[fill=white] circle(\vr);
\draw(x24)[fill=white] circle(\vr);
\draw(x25)[fill=white] circle(\vr);
\node at (5.5,0.2) {$e$};
\node at (4,-0.7) {$f$};
\node at (7,-0.7) {$f'$};
\end{scope}
\end{tikzpicture}
\caption{Graph $G_{k}'$.}
	\label{fig:graph G'}
\end{center}
\end{figure}

Let $e$ be the edge incident with $f$ and $f'$. Then it is straightforward to check that $\gp(G_k') = 4k$ and $\gp(G_k'-e) = 2k$.

To prove sharpness of the upper bound, let $k\geq 3$ and let $G_{k}''$ be a graph constructed similarly as $G_{k}'$, the construction of $G_{k}''$ should be clear from Fig.~\ref{fig:graph G''}.

\begin{figure}[ht!]
\begin{center}
\begin{tikzpicture}[scale=1.0,style=thick,x=1cm,y=1cm]
\def\vr{2pt}

\begin{scope}[xshift=0cm, yshift=0cm] 
\coordinate(x1) at (0,0);
\coordinate(x2) at (1,0.4);
\coordinate(x3) at (2,0);
\coordinate(x4) at (3,0);
\coordinate(x5) at (5,0);
\coordinate(x6) at (6,0);
\coordinate(x7) at (8,0);
\coordinate(x8) at (9,0);
\coordinate(x9) at (10,0.4);
\coordinate(x10) at (11,0);
\coordinate(x11) at (11,2);
\coordinate(x12) at (8.5,2);
\coordinate(x13) at (7,2);
\coordinate(x14) at (5.5,2);
\coordinate(x15) at (4,2);
\coordinate(x16) at (2.5,2);
\coordinate(x17) at (0,2);
\coordinate(x18) at (1,0.75);
\coordinate(x19) at (1,1.5);
\coordinate(x20) at (4,0.75);
\coordinate(x21) at (4,1.5);
\coordinate(x22) at (7,0.75);
\coordinate(x23) at (7,1.5);
\coordinate(x24) at (10,0.75);
\coordinate(x25) at (10,1.5);
\coordinate(w1) at (4,0.4);
\coordinate(w2) at (7,0.4);

\draw (x1) -- (x2) -- (x3) -- (x4);
\draw (x5)-- (x6);
\draw (x7)--(x8) -- (x9) -- (x10) -- (x11) --(x12)-- (x13) --(x14)--(x15) -- (x16) -- (x17) -- (x1);
\draw (x1) -- (x18) -- (x3) -- (x19)--(x1);
\draw (x4) -- (x20) -- (x5) -- (x21)--(x4);
\draw (x6) -- (x22) -- (x7) -- (x23)--(x6);
\draw (x8) -- (x24) -- (x10) -- (x25)--(x8);
\draw (1,1.2) node{\vdots};
\draw (4,1.2) node{\vdots};
\draw (7,1.2) node{\vdots};
\draw (10,1.2) node{\vdots};
\draw (x4) -- (w1) -- (x5);
\draw (x6) -- (w2) -- (x7);
\draw (x1) .. controls (0.5,-0.5) and (1.5,-0.5) .. (x3);
\draw (x8) .. controls (9.5,-0.5) and (10.5,-0.5) .. (x10);

\draw(w1)[fill=white] circle(\vr);
\draw(w2)[fill=white] circle(\vr);
\draw(x1)[fill=white] circle(\vr);
\draw(x2)[fill=white] circle(\vr);
\draw(x3)[fill=white] circle(\vr);
\draw(x4)[fill=white] circle(\vr);
\draw(x5)[fill=white] circle(\vr);
\draw(x6)[fill=white] circle(\vr);
\draw(x7)[fill=white] circle(\vr);
\draw(x8)[fill=white] circle(\vr);
\draw(x9)[fill=white] circle(\vr);
\draw(x10)[fill=white] circle(\vr);
\draw(x11)[fill=white] circle(\vr);
\draw(x12)[fill=white] circle(\vr);
\draw(x13)[fill=white] circle(\vr);
\draw(x14)[fill=white] circle(\vr);
\draw(x15)[fill=white] circle(\vr);
\draw(x16)[fill=white] circle(\vr);
\draw(x17)[fill=white] circle(\vr);
\draw(x18)[fill=white] circle(\vr);
\draw(x19)[fill=white] circle(\vr);
\draw(x20)[fill=white] circle(\vr);
\draw(x21)[fill=white] circle(\vr);
\draw(x22)[fill=white] circle(\vr);
\draw(x23)[fill=white] circle(\vr);
\draw(x24)[fill=white] circle(\vr);
\draw(x25)[fill=white] circle(\vr);
\node at (5.5,0.2) {$e$};
\end{scope}

\end{tikzpicture}
\caption{Graph $G_{k}''$.}
	\label{fig:graph G''}
\end{center}
\end{figure}

We can directly check that $\gp(G_k'') = 2k$ and $\gp(G_k''-e) = 4k$.
\qed

We next consider how the general position number changes when removing an edge from diameter $2$ graphs.

\begin{theorem}
\label{thm:edge-removal diameter 2}
If $e$ is an edge of a diameter $2$ graph $G$, then
$$ \gp(G)-1\le \gp(G-e)\leq\ \gp(G) + 1\,.$$
Moreover, the bounds are sharp.
\end{theorem}

\proof
Let $e=uv$ and let $X$ be an arbitrary $\gp$-set of $G$. Assume first that $u,v\notin X$. Then we claim that $X$ is also a general position set in $G-e$. To see it, consider any two vertices $x, y\in X$. Since $\diam(G) = 2$ and $u,v\notin X$, no shortest $x,y$-path of $G$ contains the edge $uv$, hence ${\cal P}_{G}(x,y) = {\cal P}_{G-e}(x,y)$. Thus $X$ is a general position set of $G-e$, so that in this case $\gp(G-e) \ge \gp(G)$. Assume second that $u\in X$. Then we claim that $X' = X\setminus \{u\}$ is a general position set in $G-e$. If $v\notin X$,  then we see by the same argument as above that $X'$ is a general position set of $G-e$. It remains to consider the subcase when $v\in X$. Since $\diam(G) = 2$, the only way that $X'$ is not a general position set in $G-e$ would be when there is a shortest path $Q$ in $G-e$ containing $v$ and two other vertices of $X'$. Since the path $Q$ does not contain the vertex $u$, $Q$ is also a shortest path in $G$, impossible. We can conclude that $X'$ is a general position set of $G-e$ and therefore $\gp(G-e) \ge \gp(G) - 1$.


Let $Y$ be an arbitrary $\gp$-set of $G-e$. Assume that $u,v\notin Y$. We claim that $Y$ is a general position set of $G$.
Consider any two vertices $x, y$ from $Y$ and let $P$ be an arbitrary shortest $x,y$-path in $G$.
Then the path $P$ does not contain the edge $uv$ in $G$ as $\diam(G) = 2$, hence ${\cal P}_{G}(x,y) = {\cal P}_{G-e}(x,y)$. We thus have that $Y$ is a general position set of $G$, and $\gp(G)\geq \gp(G-e)$.
Assume now that $u\in Y$. Then we claim that $Y' = Y\setminus \{u\}$ is a general position set of $G$. If $v\notin Y$, then similarly as above, $Y'$ is a general position set of $G$. It remains to consider the case when $v\in Y$. Suppose on the contrary that there are three vertices from $Y'$ such that they lie on a common shortest path $P'$ in $G$. Then the path $P'$ contains $v$ and two other vertices from $Y'$. The path $P'$ does not contain the edge $e$ and so $P'$ is also a shortest path in $G-e$, a contradiction. Hence $Y'$ is a general position set of $G$ and we can conclude that $\gp(G)\geq \gp(G-e) - 1$.

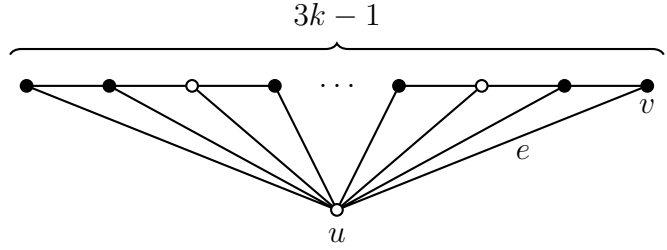
\begin{figure}[ht!]
\begin{center}
\begin{tikzpicture}[scale=1.1,style=thick,x=1cm,y=1cm]
\def\vr{2pt}

\begin{scope}[xshift=0cm, yshift=0cm] 
\coordinate(x1) at (0.25,0);
\coordinate(x2) at (-3.5,1.5);
\coordinate(x3) at (-2.5,1.5);
\coordinate(x4) at (-1.5,1.5);
\coordinate(x5) at (-0.5,1.5);
\coordinate(x6) at (1,1.5);
\coordinate(x7) at (2,1.5);
\coordinate(x8) at (3,1.5);
\coordinate(x9) at (4,1.5);

\draw (x2) --(x3)-- (x4)--(x5);
\draw (x6)--(x7) -- (x8)--(x9);
\draw (x2) -- (x1) -- (x3);
\draw (x4) -- (x1) -- (x5);
\draw (x6) -- (x1) -- (x7);
\draw (x8) -- (x1) -- (x9);
\draw (0.25,1.5) node{$\ldots$};
\draw [decorate, decoration = {brace, raise=4pt, amplitude=0.15cm}, black!100] (-3.7,1.75)--(4.2,1.75);

\draw(x1)[fill=white] circle(\vr);
\draw(x2)[fill=black] circle(\vr);
\draw(x3)[fill=black] circle(\vr);
\draw(x4)[fill=white] circle(\vr);
\draw(x5)[fill=black] circle(\vr);
\draw(x6)[fill=black] circle(\vr);
\draw(x7)[fill=white] circle(\vr);
\draw(x8)[fill=black] circle(\vr);
\draw(x9)[fill=black] circle(\vr);

\node at (0.25,-0.3) {$u$};
\node at (4,1.25) {$v$};
\node at (2.5,0.7) {$e$};
\node at (0.25,2.35) {$3k-1$};

\end{scope}

\end{tikzpicture}
\caption{Fan graph $F_{n}$ with $n=3k-1$.}
\label{fig:fan graph}
\end{center}
\end{figure}

To show that the lower bound is sharp, consider again fan graphs $F_n$, where $n = 3k-1$, $k\geq 3$, see Fig.~\ref{fig:fan graph}.

By Proposition~\ref{prop:fan} we known that $\gp(F_{3k-1}) = \lceil \frac{2(3k-1+1)}{3}\rceil = 2k$, see~Fig.~\ref{fig:fan graph} where the black vertices form a $\gp$-set of $F_n$. Let $uv$ be the edge of $F_n$ as shown in the figure. Then it is straightforward to check that all the black vertices but $v$ form a $\gp$-set of $F_n-e$, hence $\gp(F_n-e) = \gp(F_n) -1 = 2k-1$.

To show sharpness of the upper bound, let $m\ge 3$ and let $G_m$ be the graph constructed as follows. Start with $K_m$ and let $x,y$ be arbitrary, fixed vertices of it. Then we set:
\begin{align*}
V(G_m) & = V(K_m) \cup \{x',y'\}\,,\\
E(G_m) & = E(K_m) \cup \{x'y', x'x, y'y, x'y\}\,.
\end{align*}
Set $e=xx'$. Then $\diam(G_m) = 2$ and it is easy to verify that $\gp(G_m-e) = m + 1 = \gp(G_m) + 1$.
\qed

\section{Concluding remarks}

In this paper we have explored how much the general position number of a graph can be affected by removing a vertex or by removing an edge. There are several problems that seem interesting to investigate in the future. 

First, we have proved several lower and upper bounds and demonstrated that they are sharp. A challenging task remaining is to characterize the equality cases. It would also be of interest to determine sharper bounds for higher diameters. Two additional fundamental local graph operations which could be worth investigation with respect to the general position number are edge contraction and edge subdivision. Finally, the same line of research as in this paper could be done for different versions of the general position problem that appear in the literature, such as monophonic position sets. 

\section*{Acknowledgements}

We would like to thank one of the reviewers for her/his accurate reading and many helpful tips. Among other things, the advice led to a shorter proof of Theorem~\ref{thm:upper bound}.

This work was supported by the Slovenian Research Agency ARIS (research core funding P1-0297 and projects N1-0285, N1-0355).

\section*{Declaration of interests}

The authors declare that they have no known competing financial interests or personal relationships that could have appeared to influence the work reported in this paper.

\section*{Data availability}
Our manuscript has no associated data.


\end{document}